\documentclass[12pt,a4paper,onecolumn]{IEEEtran}

\usepackage{comment,url}
\usepackage[dvips]{graphicx}
\usepackage{caption,subfig}
\usepackage{float}
\graphicspath{{./fig/}}
\DeclareGraphicsExtensions{.eps}
\usepackage{setspace}
\usepackage{amsmath}
\usepackage{amssymb}
\usepackage{algorithm}
\usepackage{algorithmic}
\usepackage{cite}
\usepackage{array}
\usepackage{mdwmath}
\usepackage{mdwtab}

\allowdisplaybreaks[4]
\long\def\symbolfootnote[#1]#2{\begingroup%
\def\thefootnote{\fnsymbol{footnote}}\footnote[#1]{#2}\endgroup}

\newcommand{\dd}{\mathrm{d}}
\newcommand{\n}{\diamond}
\newcommand{\A}{\tilde{A}}
\begin{document}

\title{\LARGE{Random Distances Associated with Arbitrary Triangles: A Systematic Approach between Two Random Points}}

\author{Fei Tong, Maryam Ahmadi, and Jianping Pan\\
University of Victoria, Victoria, BC, Canada}

\maketitle

\begin{abstract}
It has been known that the distribution of the random distances between two uniformly distributed points within a convex polygon can be obtained based on its chord length distribution (CLD). In this report, we first verify the existing known CLD for arbitrary triangles, and then derive and verify the distance distribution between two uniformly distributed points within an arbitrary triangle by simulation. Furthermore, a decomposition and recursion approach is applied to obtain the random point distance distribution between two arbitrary triangles sharing a side. As a case study, the explicit distribution functions are derived when two congruent isosceles triangles with the acute angle equal to $\frac{\pi}{6}$ form a rhombus or a concave 4-gon.
\end{abstract}


\section{Problem Statement}
The goal is to obtain the probability density function (PDF) and  cumulative distribution function (CDF) of the random distances between two uniformly distributed points within an arbitrary triangle, and those between two uniformly distributed points separately inside two adjacent arbitrary triangles.
\section{Random Distances within an Arbitrary Triangle}\label{sec:cld_triangle}
In this section, we first introduce the chord length distribution (CLD) function for an arbitrary triangle, based on which the distance distribution between two random points within the triangle is then derived using a systematic approach. The obtained results are verified by simulation as well as by comparison with existing results.
\subsection{CLD for arbitrary triangles}\label{subsec:cld}
\begin{figure}[ht]
\centering
\includegraphics[width=2.5in]{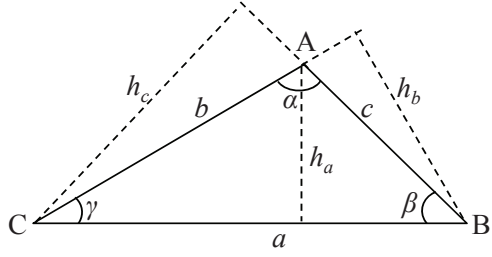}
\caption{An arbitrary triangle.}
\label{fig:triangle}
\end{figure}
An approach based on analytic geometry is applied in \cite{CLD} to obtain the CLD for arbitrary triangles. Let $\triangle{}ABC$ denote an arbitrary triangle with the side lengths $a$, $b$, and $c$, the internal angles $\alpha$, $\beta$, and $\gamma$, and the altitudes $h_a$, $h_b$, and $h_c$, respectively (see an example shown in \figurename~\ref{fig:triangle}). Without loss of generality, let $a\ge{}b\ge{}c$. The perimeter of the triangle is $u=a+b+c$, and area is $\tilde{A}=\sqrt{\frac{u}{2}(\frac{u}{2}-a)(\frac{u}{2}-b)(\frac{u}{2}-c)}$.

Let $l$ represent the chord length. According to \cite{CLD}, the CDF of the chord length of the triangle $\triangle{}ABC$ is given in three cases according to the variation of $\alpha$ and/or the relationship between $h_c$ and $c$:
\begin{figure}[H]
\begin{equation}\label{eq:cld}
  F(l)=\left\{
  \begin{array}{ll}
	F_1(l) & \alpha>\frac{\pi}{2}\\
	F_2(l) & \alpha\leq \frac{\pi}{2}~\mbox{and}~h_c<c\\
	F_3(l) & \alpha\leq \frac{\pi}{2}~\mbox{and}~h_c\geq c
  \end{array}
\right..
\end{equation}Specifically,
\end{figure}
\begin{enumerate}
  \item $\alpha>\frac{\pi}{2}$:
  \begin{figure}[H]
  \begin{equation*}
  F_1(l)=\left\{
  \begin{array}{ll}
	0 & l<0\\
	\frac{H_1(l)}{u} & 0\leq l\leq{}h_a\\
	\frac{H_2(l)}{u} & h_a\leq l\leq c\\
	\frac{H_3(l)}{u} & c\leq l\leq b\\
	\frac{H_4(l)}{u} & b\leq l\leq a\\
	1 & l>a
  \end{array}
\right.,
\end{equation*}
\end{figure}
  \item $\alpha\leq\frac{\pi}{2}$ and $h_c<c$:
  \begin{figure}[H]
  \begin{equation*}
  F_{2}(l)=\left\{
  \begin{array}{ll}
	0 & l<0\\
	\frac{H_1(l)}{u} & 0\leq l\leq{}h_a\\
	\frac{H_2(l)}{u} & h_a\leq l\leq h_b\\
	\frac{H_5(l)}{u} & h_b\leq l\leq h_c\\
    \frac{H_6(l)}{u} & h_c\leq l\leq c\\
    \frac{H_7(l)}{u} & c\leq l\leq b\\
	\frac{H_4(l)}{u} & b\leq l\leq a\\
	1 & l>a
  \end{array}
\right.,
\end{equation*}
\end{figure}
  \item $\alpha\leq\frac{\pi}{2}$ and $h_c\geq c$:
  \begin{figure}[H]
  \begin{equation*}
    F_3(l)=\left\{
    \begin{array}{ll}
	   0 & l<0\\
       \frac{H_1(l)}{u} & 0\leq l\leq{}h_a\\
	   \frac{H_2(l)}{u} & h_a\leq l\leq h_b\\
       \frac{H_5(l)}{u} & h_b\leq l\leq c\\
	   \frac{H_3(l)}{u} & c\leq l\leq h_c\\
       \frac{H_7(l)}{u} & h_c\leq l\leq b\\
	   \frac{H_4(l)}{u} & b\leq l\leq a\\
	   1 & l>a
    \end{array}
    \right.,
 \end{equation*}\end{figure}
\end{enumerate}%
where
\begin{align*}
  H_1(l) = {} & \frac{3l}{2}+\frac{l}{2}[(\pi-\alpha)\cot\alpha+(\pi-\beta)\cot\beta+(\pi-\gamma)\cot\gamma] \,,\\
  H_2(l) = {} & \frac{3l}{2}+a\sin\varphi_1+\frac{l}{2}\left[(\pi-\alpha)\cot\alpha+(\pi-\beta-2\varphi_1)\cot\beta+(\pi-\gamma-2\varphi_1)\cot\gamma\right] \,,\\
  H_3(l) = {} & l+c+\frac{a}{2}\sin\varphi_1+\frac{b}{2}\sin\varphi_2+\frac{l}{2}\left[\left(\frac{\pi}{2}-\varphi_2\right)\cot\alpha\right.\\
  &\left.+\left(\frac{\pi}{2}-\varphi_1\right)\cot\beta+(\pi-2\gamma-\varphi_1-\varphi_2)\cot\gamma\right] \,,\\
  H_4(l) = {} & \frac{l}{2}+b+c+\frac{b}{2}\sin\varphi_2+\frac{c}{2}\sin\varphi_3+\frac{l}{2}\left[(\alpha-\varphi_2-\varphi_3)\cot\alpha+\left(\frac{\pi}{2}-\beta-\varphi_3\right)\cot\beta\right.\\
  &\left.+\left(\frac{\pi}{2}-\gamma-\varphi_2\right)\cot\gamma\right] \,,\\
  H_5(l) = {} & \frac{3l}{2}+a\sin\varphi_1+b\sin\varphi_2+\frac{l}{2}\left[(\pi-\alpha-2\varphi_2)\cot\alpha+(\pi-\beta-2\varphi_1)\cot\beta\right.\\
  &\left.+(\pi-\gamma-2\varphi_1-2\varphi_2)\cot\gamma\right] \,,\\
  H_6(l) = {} & \frac{3l}{2}+a\sin\varphi_1+b\sin\varphi_2+c\sin\varphi_3+\frac{l}{2}\left[(\pi-\alpha-2\varphi_2-2\varphi_3)\cot\alpha\right.\\
  &\left.+(\pi-\beta-2\varphi_1-2\varphi_3)\cot\beta+(\pi-\gamma-2\varphi_1-2\varphi_2)\cot\gamma\right] \,,\\
  H_7(l) = {} & l+c+\frac{a}{2}\sin\varphi_1+\frac{b}{2}\sin\varphi_2+c\sin\varphi_3+\frac{l}{2}\left[\left(\frac{\pi}{2}-\varphi_2-2\varphi_3\right)\cot\alpha\right.\\
  &+\left.\left(\frac{\pi}{2}-\varphi_1-2\varphi_3\right)\cot\beta+(\pi-2\gamma-\varphi_1-\varphi_2)\cot\gamma\right] \,,
\end{align*}
with $\varphi_1=\arccos\frac{h_a}{l}$, $\varphi_2=\arccos\frac{h_b}{l}$, and $\varphi_3=\arccos\frac{h_c}{l}$.

In \cite{CLD}, the above chord length CDF in (\ref{eq:cld}) was not verified, which will be done in Section~\ref{subsec:verify} through simulation.

\subsection{Point distance distribution for arbitrary triangles}\label{subsec:pdd_dist}
In this subsection, we derive the PDF and CDF of the distances between two random points within an arbitrary triangle.
\subsubsection{PDF}Denote the PDF of the chord length $l$ for $\triangle{}ABC$ by $f(l)$, and the PDF of the distance $d~(0\leq d\leq a)$ between two random points in $\triangle{}ABC$ by $g(d)$. According to \cite{piefke1978beziehungen}, the relationship between these two functions is given by
\begin{equation}\label{eq:g}
  g(d) = \frac{2ud}{\tilde{A}^2}\int_d^a (l-d)f(l)\,\mathrm{d} l\,.
\end{equation}
Following \cite{basel2012random}, we have
\begin{equation*}
 \begin{array}{ll}
	g(d){}& =~ \frac{2d}{\tilde{A}}\left[\pi-\frac{u}{\tilde{A}}\left(d-\int_0^dF(l)\,\mathrm{d} l\right)\right]\\
    	&=~\frac{2d}{\tilde{A}}\left[\pi+\frac{1}{\tilde{A}}\,(I^*(d)-ud)\right]\,,
    \end{array}
\end{equation*}%
where $I^*(d)=u\int_0^d F(l)\,\mathrm{d} l$. With (\ref{eq:cld}) we have
\begin{figure}[H]
\begin{equation}\label{eq:pdf_Istar}
  I^*(d)=\left\{
  \begin{array}{ll}
	I^*_1(d) & \alpha>\frac{\pi}{2}\\
	I^*_2(d) & \alpha\leq \frac{\pi}{2}~\mbox{and}~h_c<c\\
	I^*_3(d) & \alpha\leq \frac{\pi}{2}~\mbox{and}~h_c\geq c
  \end{array}
\right..
\end{equation}
\end{figure}
\noindent{}Specifically,
\begin{enumerate}
  \item $\alpha>\frac{\pi}{2}$:
  \begin{figure}[H]
  \begin{equation*}
  I^*_1(d)=\left\{
  \begin{array}{ll}
	J_1^*(0,d) & 0\leq d\leq{}h_a\\
	J_1^*(0,h_a)+J_2^*(h_a,d) & h_a\leq d\leq c\\
	J_1^*(0,h_a)+J_2^*(h_a,c)+J_3^*(c,d) & c\leq d\leq b\\
	J_1^*(0,h_a)+J_2^*(h_a,c)+J_3^*(c,b)+ J_4^*(b,d)& b\leq d\leq a
  \end{array}
\right.,
\end{equation*}
\end{figure}
  \item $\alpha\leq\frac{\pi}{2}$ and $h_c<c$:
  \begin{figure}[H]
  \begin{equation*}
  I^*_{2}(d)=\left\{
  \begin{array}{ll}
	J_1^*(0,d) & 0\leq d\leq{}h_a\\
	J_1^*(0,h_a)+J_2^*(h_a,d) & h_a\leq d\leq h_b\\
	J_1^*(0,h_a)+J_2^*(h_a,h_b)+J_5^*(h_b,d) & h_b\leq d\leq h_c\\
    J_1^*(0,h_a)+J_2^*(h_a,h_b)+J_5^*(h_b,h_c)+J_6^*(h_c,d) & h_c\leq d\leq c\\
    J_1^*(0,h_a)+J_2^*(h_a,h_b)+J_5^*(h_b,h_c)+J_6^*(h_c,c)+J_7^*(c,d) & c\leq d\leq b\\
	J_1^*(0,h_a)+J_2^*(h_a,h_b)+J_5^*(h_b,h_c)+J_6^*(h_c,c)+J_7^*(c,b)\\
~~~+J_4^*(b,d) & b\leq d\leq a
  \end{array}
\right.,
\end{equation*}
\end{figure}
  \item $\alpha\leq\frac{\pi}{2}$ and $h_c\geq c$:
  \begin{figure}[H]\begin{equation*}
    I^*_3(d)=\left\{
    \begin{array}{ll}
       J_1^*(0,d) & 0\leq d\leq{}h_a\\
	   J_1^*(0,h_a)+J_2^*(h_a,d) & h_a\leq d\leq h_b\\
       J_1^*(0,h_a)+J_2^*(h_a,h_b)+J_5^*(h_b,d) & h_b\leq d\leq c\\
	   J_1^*(0,h_a)+J_2^*(h_a,h_b)+J_5^*(h_b,c)+J_3^*(c,d) & c\leq d\leq h_c\\
       J_1^*(0,h_a)+J_2^*(h_a,h_b)+J_5^*(h_b,c)+J_3^*(c,h_c)+J_7^*(h_c,d) & h_c\leq d\leq b\\
	   J_1^*(0,h_a)+J_2^*(h_a,h_b)+J_5^*(h_b,c)+J_3^*(c,h_c)+J_7^*(h_c,b)\\
~~~+J_4^*(b,d) & b\leq d\leq a
    \end{array}
    \right.,
 \end{equation*}
\end{figure}
\end{enumerate}%
with
\begin{equation}\label{eq:J_star}
  J_k^*(l,d) = H_k^*(d) - H_k^*(l)\,,k=1\dots{}7,
\end{equation}%
\begin{equation}\label{eq:H_star}
  H_k^*(d) = \int{}H_k(d)\dd d\,,k=1\dots{}7,
\end{equation}
where
\begin{align*}
  H_1^*(d) = {} & \frac{{d}^{2}}{4} \left[ (\pi-\alpha)\cot\alpha   + (\pi-\beta)\cot\beta  + (\pi-\gamma)\cot\gamma+3\right]
 \,,\\
  H_2^*(d) = {} & \frac{3a}{2}\sqrt{d^2-h_a^2}-\frac{ad^2}{2h_a}\arccos\frac{h_a}{d}+ah_a\arcsin\frac{h_a}{d}\\
  &+\frac{{d}^{2}}{4} \left[ (\pi-\alpha)\cot\alpha   + (\pi-\beta)\cot\beta  + (\pi-\gamma)\cot\gamma+3\right]\,,\\
  H_3^*(d) = {} &  \frac{3}{4}\left(a\sqrt{d^2-h_a^2}+b\sqrt{d^2-h_b^2}\right)-\frac{d^2}{4}\left(\frac{a}{h_a}\arccos\frac{h_a}{d}+\frac{b}{h_b}\arccos\frac{h_b}{d}\right)\\
  &+\frac{1}{2}\left( ah_a\arcsin\frac{h_a}{d}+bh_b\arcsin\frac{h_b}{d} \right)+\frac{d}{2}\left[\frac{\pi{}d}{4}\left(\frac{a}{h_a}+\frac{b}{h_b}\right)+d+2c-\gamma{}d\cot\gamma\right]\,,\\
  H_4^*(d) = {} & \frac{3}{4}\left(b\sqrt{d^2-h_b^2}+c\sqrt{d^2-h_c^2}\right) - \frac{d^2}{4}\left(\frac{b}{h_b}\arccos\frac{h_b}{d}+\frac{c}{h_c}\arccos\frac{h_c}{d}\right)\\
  &+\frac{1}{2}\left( bh_b\arcsin\frac{h_b}{d}+ch_c\arcsin\frac{h_c}{d} \right)\\
  &+d\left[\frac{d}{4}(\alpha\cot\alpha-\beta\cot\beta-\gamma\cot\gamma)+\frac{\pi{}ad}{8h_a}+\frac{d}{4}+b+c\right]\,,\\
  H_5^*(d) = {} & \frac{3}{2}\left(a\sqrt{d^2-h_a^2}+b\sqrt{d^2-h_b^2}\right)- \frac{d^2}{2}\left(\frac{a}{h_a}\arccos\frac{h_a}{d}+\frac{b}{h_b}\arccos\frac{h_b}{d}\right)\\ &+ah_a\arcsin\frac{h_a}{d}+bh_b\arcsin\frac{h_b}{d}\\
  &+\frac{{d}^{2}}{4} \left[ (\pi-\alpha)\cot\alpha   + (\pi-\beta)\cot\beta  + (\pi-\gamma)\cot\gamma+3\right]\,,\\
  H_6^*(d) = {} & \frac{3}{2}\left(a\sqrt{d^2-h_a^2}+b\sqrt{d^2-h_b^2}+c\sqrt{d^2-h_c^2}\right)\\
  &-\frac{d^2}{2}\left(\frac{a}{h_a}\arccos\frac{h_a}{d}+\frac{b}{h_b}\arccos\frac{h_b}{d}+\frac{c}{h_c}\arccos\frac{h_c}{d}\right)\\
  &+ah_a\arcsin\frac{h_a}{d}+bh_b\arcsin\frac{h_b}{d}+ch_c\arcsin\frac{h_c}{d}\\
  &+\frac{{d}^{2}}{4} \left[ (\pi-\alpha)\cot\alpha   + (\pi-\beta)\cot\beta  + (\pi-\gamma)\cot\gamma+3\right]\,,\\
  H_7^*(d) = {} & \frac{3}{4}\left(a\sqrt{d^2-h_a^2}+b\sqrt{d^2-h_b^2}+2c\sqrt{d^2-h_c^2}\right)\\
  &-\frac{d^2}{4}\left(\frac{a}{h_a}\arccos\frac{h_a}{d}+\frac{b}{h_b}\arccos\frac{h_b}{d}+\frac{2c}{h_c}\arccos\frac{h_c}{d}\right)\\
  &+\frac{1}{2}\left( ah_a\arcsin\frac{h_a}{d}+bh_b\arcsin\frac{h_b}{d}+2ch_c\arcsin\frac{h_c}{d}\right)\\
  &+\frac{d}{2}\left[ \frac{\pi{}d}{4}\left(\frac{a}{h_a}+\frac{b}{h_b}\right)+d+2c-\gamma{}d\cot\gamma\right]\,.
\end{align*}
\subsubsection{CDF}The CDF $G$ of the distance $d~(0\leq{}d\leq{}a)$ between two random points in $\triangle{}ABC$ is
\begin{align}\label{eq:G}
\begin{aligned}
  G(d)
	= {} & \int_0^d g(\tau)\,\mathrm{d}\tau
	= \int_0^d\left(\frac{2\pi\tau}{\A}-\frac{2u\tau^2}{\A^2}+\frac{2u\tau}{\A^2}\int_0^\tau F(l)\,\dd l\right)\dd\tau\\
	= {} & \frac{\pi d^2}{\A}-\frac{2ud^3}{3\A^2}+\frac{2}{\A^2}\int_0^d\tau\left(u\int_0^\tau F(l)\,\dd l\right)\dd\tau\\
	= {} & \frac{\pi d^2}{\A}-\frac{2ud^3}{3\A^2}+\frac{2}{\A^2}\int_0^d\tau I^*(\tau)\,\dd\tau\\
	= {} & \frac{1}{\A}\left[d^2\left(\pi-\dfrac{2u}{3\A}\,d\right)+\dfrac{2}{\A}\,I^\n(d)\right]\,,
\end{aligned}
\end{align}%
where $I^\n(d)=\int_0^d\tau I^*(\tau)\,\dd\tau$. With (\ref{eq:pdf_Istar}), we have
\begin{figure}[H]
\begin{equation}\label{eq:I_diamond}
  I^\n(d)=\left\{
  \begin{array}{ll}
	I^\n_1(d) & \alpha>\frac{\pi}{2}\\
	I^\n_2(d) & \alpha\leq \frac{\pi}{2}~\mbox{and}~h_c<c\\
	I^\n_3(d) & \alpha\leq \frac{\pi}{2}~\mbox{and}~h_c\geq c
  \end{array}
\right..
\end{equation}
\end{figure}
\noindent{}Specifically,
\begin{enumerate}
  \item $\alpha>\frac{\pi}{2}$:
\begin{figure}[H]
  \begin{equation*}
  I^\n_1(d)=\left\{
  \begin{array}{ll}
	K_{11}(d) & 0\leq d\leq{}h_a\\
	K_{11}(h_a)+K_{12}(d) & h_a\leq d\leq c\\
	K_{11}(h_a)+K_{12}(c)+K_{13}(d) & c\leq d\leq b\\
	K_{11}(h_a)+K_{12}(c)+K_{13}(b)+ K_{14}(d)& b\leq d\leq a
  \end{array}
\right.,
\end{equation*}
\end{figure} where $K_{ik}$ indicates that it is for case ``$i$)'' ($i=1,~2,~\text{or}~3$; e.g., if $i=1$, it is the current case $\alpha>\frac{\pi}{2}$), and $k=1\dots7$:
\begin{align*}
  K_{11}(d) = {} & J_1^\n(0,d)\,,\\
  K_{12}(d) = {} & \frac{1}{2}(d^2-h_a^2)\left[J_1^*(0,h_a)-H_2^*(h_a)\right]+J_2^\n(h_a,d)\,,\\
  K_{13}(d) = {} & \frac{1}{2}(d^2-c^2)\left[J_1^*(0,h_a)+J_2^*(h_a,c)-H_3^*(c)\right]+J_3^\n(c,d)\,,\\
  K_{14}(d) = {} & \frac{1}{2}(d^2-b^2)\left[J_1^*(0,h_a)+J_2^*(h_a,c)+J_3^*(c,b)-H_4^*(b)\right]+J_4^\n(b,d)\,.
\end{align*}
  \item $\alpha\leq\frac{\pi}{2}$ and $h_c<c$:
  \begin{figure}[H]
  \begin{equation*}
  I^\n_{2}(d)=\left\{
  \begin{array}{ll}
	K_{21}(d) & 0\leq d\leq{}h_a\\
	K_{21}(h_a)+K_{22}(d) & h_a\leq d\leq h_b\\
	K_{21}(h_a)+K_{22}(h_b)+K_{25}(d) & h_b\leq d\leq h_c\\
    K_{21}(h_a)+K_{22}(h_b)+K_{25}(h_c)+K_{26}(d) & h_c\leq d\leq c\\
    K_{21}(h_a)+K_{22}(h_b)+K_{25}(h_c)+K_{26}(c)+K_{27}(d) & c\leq d\leq b\\
	K_{21}(h_a)+K_{22}(h_b)+K_{25}(h_c)+K_{26}(c)+K_{27}(b)+K_{24}(d) & b\leq d\leq a
  \end{array}
\right.,
\end{equation*}
\end{figure}where
\begin{align*}
  K_{21}(d) = {} & K_{11}(d)=J_1^\n(0,d)\,,\\
  K_{22}(d) = {} & K_{12}(d)=\frac{1}{2}(d^2-h_a^2)\left[J_1^*(0,h_a)-H_2^*(h_a)\right]+J_2^\n(h_a,d)\,,\\
  K_{25}(d) = {} & \frac{1}{2}(d^2-h_b^2)\left[J_1^*(0,h_a)+J_2^*(h_a,h_b)-H_5^*(h_b)\right]+J_5^\n(h_b,d)\,,\\
  K_{26}(d) = {} & \frac{1}{2}(d^2-h_c^2)\left[J_1^*(0,h_a)+J_2^*(h_a,h_b)+J_5^*(h_b,h_c)-H_6^*(h_c)\right]+J_6^\n(h_c,d)\,,\\
  K_{27}(d) = {} & \frac{1}{2}(d^2-c^2)\left[J_1^*(0,h_a)+J_2^*(h_a,h_b)+J_5^*(h_b,h_c)+J_6^*(h_c,c)-H_7^*(c)\right]+J_7^\n(c,d)\,,\\
  K_{24}(d) = {} & \frac{1}{2}(d^2-b^2)\left[J_1^*(0,h_a)+J_2^*(h_a,h_b)+J_5^*(h_b,h_c)+J_6^*(h_c,c)+J_7^*(c,b)-H_4^*(b)\right]+J_4^\n(b,d)\,.
\end{align*}

  \item $\alpha\leq\frac{\pi}{2}$ and $h_c\geq c$:
  \begin{figure}[H]\begin{equation*}
    I^\n_3(d)=\left\{
    \begin{array}{ll}
       K_{31}(d) & 0\leq d\leq{}h_a\\
	   K_{31}(h_a)+K_{32}(d) & h_a\leq d\leq h_b\\
       K_{31}(h_a)+K_{32}(h_b)+K_{35}(d) & h_b\leq d\leq c\\
	   K_{31}(h_a)+K_{32}(h_b)+K_{35}(c)+K_{33}(d) & c\leq d\leq h_c\\
       K_{31}(h_a)+K_{32}(h_b)+K_{35}(c)+K_{33}(h_c)+K_{37}(d) & h_c\leq d\leq b\\
	   K_{31}(h_a)+K_{32}(h_b)+K_{35}(c)+K_{33}(h_c)+K_{37}(b)+K_{34}(d) & b\leq d\leq a
    \end{array}
    \right.,
 \end{equation*}
\end{figure}where
\begin{align*}
  K_{31}(d) = {} & K_{21}(d)=K_{11}(d)= J_1^\n(0,d)\,,\\
  K_{32}(d) = {} & K_{22}(d)=K_{12}(d)=\frac{1}{2}(d^2-h_a^2)\left[J_1^*(0,h_a)-H_2^*(h_a)\right]+J_2^\n(h_a,d)\,,\\
  K_{35}(d) = {} & K_{25}(d)= \frac{1}{2}(d^2-h_b^2)\left[J_1^*(0,h_a)+J_2^*(h_a,h_b)-H_5^*(h_b)\right]+J_5^\n(h_b,d)\,,\\
  K_{33}(d) = {} & \frac{1}{2}(d^2-c^2)\left[J_1^*(0,h_a)+J_2^*(h_a,h_b)+J_5^*(h_b,c)-H_3^*(c)\right]+J_3^\n(c,d)\,,\\
  K_{37}(d) = {} & \frac{1}{2}(d^2-h_c^2)\left[J_1^*(0,h_a)+J_2^*(h_a,h_b)+J_5^*(h_b,c)+J_3^*(c,h_c)-H_7^*(h_c)\right]+J_7^\n(h_c,d)\,,\\
  K_{34}(d) = {} & \frac{1}{2}(d^2-b^2)\left[J_1^*(0,h_a)+J_2^*(h_a,h_b)+J_5^*(h_b,c)+J_3^*(c,h_c)+J_7^*(h_c,b)-H_4^*(b)\right]+J_4^\n(b,d)\,.
\end{align*}
\end{enumerate}%
$J_k^*$ and $H_k^*$ have been given in (\ref{eq:J_star}) and (\ref{eq:H_star}), respectively, and we have
\begin{equation}\label{eq:J_diamond}
  J_k^\n(l,d) = H_k^\n(d) - H_k^\n(l)\,,k=1\dots{}7,
\end{equation}
\begin{equation}\label{eq:H_diamond}
  H_k^\n(d) = \int{}dH_k^*(d)\dd d\,,k=1\dots{}7,
\end{equation}where
\begin{align*}
  H_1^\n(d) = {} & \frac{d^4}{16}\left[(\pi-\alpha)\cot\alpha+(\pi-\beta)\cot\beta+(\pi-\gamma)\cot\gamma+3\right]\,,\\
  H_2^\n(d) = {} & \frac{1}{48h_a}\left[(26ah_ad^2+4ah_a^3)\sqrt{d^2-h_a^2}-3d^2\left(2a d^2\arccos\frac{h_a}{d}\right.\right.\\
  &\left.\left.+h_a\left(d^2((\alpha-\pi)\cot\alpha+(\beta-\pi)\cot\beta+(\gamma-\pi)\cot\gamma-3)-8ah_a\arcsin\frac{h_a}{d}\right) \right)\right]\,,\\
  H_3^\n(d) 
  ={}&\frac{1}{96h_ah_b}\left[26ah_ah_b\left(d^2+\frac{2h_a^2}{13}\right)\sqrt{d^2-h_a^2}+ 26bh_ah_b\left(d^2+\frac{2h_b^2}{13}\right)\sqrt{d^2-h_b^2}\right.\\
  &\left.-12d^2\left(\frac{d^2}{2}\left(ah_b\arccos\frac{h_a}{d}+bh_a\arccos\frac{h_b}{d}\right)-2h_ah_b\left(ah_a\arcsin\frac{h_a}{d}\right.\right.\right. \\
  &\left.\left.\left.+bh_b\arcsin\frac{h_b}{d}\right)+d\left(\gamma{}h_ah_bd\cot\gamma-h_a\left(\left(\frac{8c}{3}+d\right)h_b+\frac{\pi{}bd}{4}\right)-\frac{\pi{}ah_bd}{4}\right)\right)\right]\,,\\
  H_4^\n(d) 
  ={}&\frac{1}{96h_ah_bh_c}\left[4h_ah_bh_c\left(b\left(\frac{13d^2}{2}+h_b^2\right)\sqrt{d^2-h_b^2}+c\left(\frac{13d^2}{2}+h_c^2\right)\sqrt{d^2-h_c^2}\right)\right.\\
  &\left.+6d^2\left( -d^2h_ah_cb\arccos\frac{h_b}{d}+h_b\left( -d^2ch_a\arccos\frac{h_c}{d}\right.\right.\right.\\
  &+h_c\left( 4h_a\left( bh_b\arcsin\frac{h_b}{d}+ch_c\arcsin\frac{h_c}{d}\right)\right.\\
  &+\left.\left.\left.\left.d\left(dh_a\left(\alpha\cot\alpha-\beta\cot\beta-\gamma\cot\gamma\right)+h_a\left(d+\frac{16}{3}(b+c)\right)+\frac{\pi{}ad}{2}\right)\right) \right)\right) \right]\,,\\
  H_5^\n(d) 
  ={}& \frac{1}{48h_ah_b}\left[26h_ah_b\left(a\left(d^2+\frac{2h_a^2}{13}\right)\sqrt{d^2-h_a^2}+b\left(d^2+\frac{2h_b^2}{13}\right)\sqrt{d^2-h_b^2}\right) \right.\\
  &\left.-3d^2\left( 2ah_bd^2\arccos\frac{h_a}{d}+h_a\left(2bd^2\arccos\frac{h_b}{d}-h_b\left(8ah_a\arcsin\frac{h_a}{d}\right.\right.\right.\right.\\
  &\left.\left.\left.\left.+8bh_b\arcsin\frac{h_b}{d}+d^2((\pi-\alpha)\cot\alpha+(\pi-\beta)\cot\beta+(\pi-\gamma)\cot\gamma+3) \right)\right)\right) \right]\,,\\
  H_6^\n(d) 
  ={}& \frac{1}{48h_ah_bh_c}\left[ 4h_ah_bh_c\left(\frac{13a}{2}\left(d^2+\frac{2h_a^2}{13}\right)\sqrt{d^2-h_a^2}+b\left(\frac{13d^2}{2}+h_b^2\right)\sqrt{d^2-h_b^2}\right.\right.\\
  &\left.\left.+c\left(\frac{13d^2}{2}+h_c^2\right)\sqrt{d^2-h_c^2} \right) -3d^2\left( 2ah_bh_cd^2\arccos\frac{h_a}{d}+h_a\left(2bh_cd^2\arccos\frac{h_b}{d}\right.\right.\right.\\
  &\left.\left.\left.+h_b\left( 2cd^2\arccos\frac{h_c}{d}-h_c\left( 8\left(ah_a\arcsin\frac{h_a}{d}+bh_b\arcsin\frac{h_b}{d}\right.\right.\right.\right.\right.\right.\\
  &\left.\left.\left.\left.\left.\left.+ch_c\arcsin\frac{h_c}{d}\right)+ d^2\left((\pi-\alpha)\cot\alpha+(\pi-\beta)\cot\beta+(\pi-\gamma)\cot\gamma+3\right)\right)\right)\right)\right)\right]\,,\\
  H_7^\n(d) 
  ={}& \frac{1}{96h_ah_bh_c}\left[26h_ah_bh_c\left(a\left(d^2+\frac{2h_a^2}{13}\right)\sqrt{d^2-h_a^2}+b\left(d^2+\frac{2h_b^2}{13}\right)\sqrt{d^2-h_b^2}\right.\right.\\
  &\left.\left.+2c\left(d^2+\frac{2h_c^2}{13}\right)\sqrt{d^2-h_c^2}\right)-12d^2\left( \frac{d^2}{2}\left( ah_bh_c\arccos\frac{h_a}{d}+bh_ah_c\arccos\frac{h_b}{d}\right.\right.\right.\\
  &\left.\left.\left.+2ch_ah_b\arccos\frac{h_c}{d} \right) - h_c\left( 2h_ah_b\left(ah_a\arcsin\frac{h_a}{d}+bh_b\arcsin\frac{h_b}{d}+2ch_c\arcsin\frac{h_c}{d}\right)\right.\right.\right.\\
  &\left.\left.\left.-d\left(\gamma{}h_ah_bd\cot\gamma-h_a\left( h_b\left(\frac{8c}{3}+d\right)+\frac{\pi{}db}{4}\right) - \frac{\pi{}ah_bd}{4}\right) \right)\right)\right]\,.
\end{align*}

\subsection{Verification by simulation}\label{subsec:verify}
In this subsection, we first verify the CLD function (\ref{eq:cld}) for arbitrary triangles derived in \cite{CLD}, which leads to an elementary geometry approach to obtaining CLD, and then verify our derived distribution of the distances between two random points within an arbitrary triangle, by comparing our results with those of simulation.
\subsubsection{CLD verification}
\begin{algorithm}[H]
\renewcommand{\algorithmicrequire}{\textbf{Input:}}
\renewcommand\algorithmicensure {\textbf{Output:}}
\caption{\scriptsize{Simulation Algorithm for Chord Length Computation}}\label{alg:sim}
 \begin{algorithmic}[1]
 \REQUIRE~~\\
 Parameters with regard to an arbitrary triangle $\triangle ABC$:\\
 $a, b, c, h_a, h_b, h_c, \beta, \gamma, A(x_A, y_A), B(x_B, y_B), C(x_C, y_C)$;\\
 $x_A=b\cos\gamma, y_A=h_a, x_B=a, y_B=0, x_C=0, y_C=0; $
 \ENSURE~~\\
 Chord length list $L_c$;
 \STATE $\delta\theta=\frac{\pi}{180}$; $\delta d=\frac{1}{1,000}$;
 \FOR{$\theta=0$; $\theta\leq\pi$; $\theta=\theta+\delta\theta$}
    \STATE $flag=0$; $//$ it is case 1), 3), 5) if flag=1, and case 2), 4), 6) if flag=2;
    \IF{$\theta==0$ or $\theta==\pi$}
        \STATE $flag=1$; $d=h_a$; $base=a$; $//$ case 1), \figurename~\ref{fig:cld_sim} (1);
    \ELSIF{$\theta==\gamma$}
         \STATE $flag=1$; $d=h_b$; $base=b$; $//$ case 3), \figurename~\ref{fig:cld_sim} (3);
    \ELSIF{$\theta==\pi-\beta$}
        \STATE $flag=1$; $d=h_c$; $base=c$; $//$ case 5), \figurename~\ref{fig:cld_sim} (5);
    \ELSIF{$0<\theta<\gamma$}
        \STATE $flag=2$; $//$ case 2), \figurename~\ref{fig:cld_sim} (2);
        \STATE $d_1=b\sin(\gamma-\theta)$; $d_2=a\sin\theta$; $d=d_1+d_2$;
        \STATE $y_I=\frac{a}{\cot\beta+\cot\theta}$; $x_I=y_I\cot\theta$; $//$ coordinate of point $I$;
        \STATE $base=\sqrt{(x_I-x_C)^2+(y_I-y_C)^2}$ $//$$|CI|$;
    \ELSIF{$\gamma<\theta<\pi-\beta$}
        \STATE $flag=2$; $//$ case 4), \figurename~\ref{fig:cld_sim} (4);
        \STATE $d_1=b\sin(\theta-\gamma)$; $d_2=c\sin(\theta+\beta)$; $d=d_1+d_2$;
        \STATE $x_I=\frac{b}{(\cot\gamma+\cot(\theta-\gamma))*sin\gamma}$; $y_I=0$; $//$ coordinate of point $I$;
        \STATE $base=\sqrt{(x_I-x_A)^2+(y_I-y_A)^2}$ $//$$|AI|$;
    \ELSIF{$\pi-\beta<\theta<\pi$}
        \STATE $flag=2$; $//$ case 6), \figurename~\ref{fig:cld_sim} (6);
        \STATE $d_1=a\sin(\theta)$; $d_2=-c\sin(\theta+\beta)$; $d=d_1+d_2$;
        \STATE $y_I=\frac{a}{\cot\gamma-\cot\theta}$; $x_I=y_I\cot\gamma$; $//$ coordinate of point $I$;
        \STATE $base=\sqrt{(x_I-x_B)^2+(y_I-y_B)^2}$ $//$$|BI|$;
    \ENDIF
    \IF{$flag==1$}
        \STATE $//$ case 1), 3), 5)
        \FOR{$d^\prime=0$; $d^\prime\leq d$; $d^\prime=d^\prime+\delta d$}
            \STATE $l=\frac{d^\prime\cdot base}{d}$; and insert $l$ into $L_c$;
        \ENDFOR
    \ELSIF{$flag==2$}
        \STATE $//$ case 2), 4), 6)
        \FOR{$d^\prime=0$; $d^\prime\leq d$; $d^\prime=d^\prime+\delta d$}
            \IF{$d^\prime\leq d_1$}
                \STATE $l=\frac{d^\prime\cdot base}{d_1}$; and insert $l$ into $L_c$;
            \ELSE
                \STATE $l=\frac{(d-d^\prime)\cdot base}{d_2}$; and insert $l$ into $L_c$;
            \ENDIF
        \ENDFOR
    \ENDIF
 \ENDFOR
 \end{algorithmic}
\end{algorithm}
\begin{figure}[ht]
\centering
\includegraphics[width=5.5in]{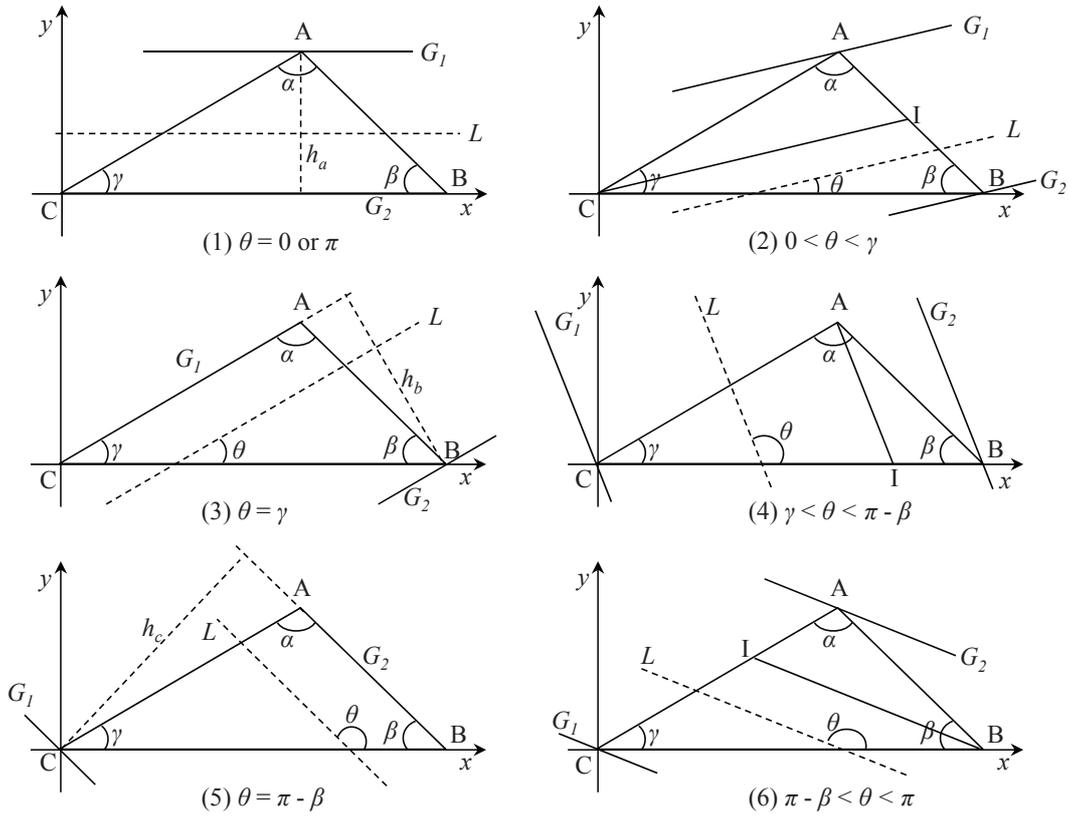}
\caption{Varying triangle chords.}
\label{fig:cld_sim}
\end{figure}

Let us build a rectangular coordinate system for the arbitrary triangle shown in \figurename~\ref{fig:triangle}. Without loss of generality, let vertex \textit{C} be located at the origin, and side $CB$ on the positive $x$-axis. For a given chord in line $L$ with a given orientation $\theta$ (without loss of generality, with regard to $CB$), there are six cases as shown in \figurename~\ref{fig:cld_sim}: 1) $\theta=0~\text{or}~\pi$; 2) $0<\theta<\gamma$; 3) $\theta=\gamma$; 4) $\gamma<\theta<\pi-\beta$; 5) $\theta=\pi-\beta$; 6) $\pi-\beta<\theta<\pi$.

The simulation is conducted as follows. Let $\theta$ increase from $0$ to $\pi$ with a fixed small step of $\delta\theta$ (e.g., $\delta\theta=\frac{\pi}{180}$). For each specific $\theta$, let the triangle be exactly between two lines $G_1$ and $G_2$ which are parallel with $L$ ($G_1$ and $G_2$ intersect the boundary of the triangle). Denote the distance between $G_1$ and $G_2$ by $d$. Using a fixed, small $\delta{}d$ (e.g., $\delta{}d=\frac{1}{1,000}$), and varying $d^\prime$ from 0 to $d$ with step $\delta d$, we can get $\frac{d}{\delta d}$ chords which are all parallel with $G_1$ and $G_2$. For each chord with the length denoted by $l$,
\begin{itemize}
  \item if it is parallel with a side of the triangle (the length of this side is denoted by $base$) as shown in \figurename~\ref{fig:cld_sim} (1), (3), and (5), $l=\frac{d^\prime\cdot base}{d}$;
  \item otherwise, $G_1$ and $G_2$ pass through two vertices of the triangle, respectively. Denote the distance from the remaining vertex to $G_1$ as $d_1$ and the distance to $G_2$ as $d_2$ ($d=d_1+d_2$). Let a line pass through the remaining vertex and be parallel with $G_1$ and $G_2$. The point at which this line and the opposite side of the remaining vertex intersects is denoted by $I$, as shown in \figurename~\ref{fig:cld_sim} (2), (4), and (6), with the corresponding cord length denoted by $base$. Without loss of generality, if $d^\prime\leq d_1$, $l=\frac{d^\prime\cdot base}{d_1}$; otherwise, $l=\frac{(d-d^\prime)\cdot base}{d_2}$.
\end{itemize}
The algorithm is summarized in Algorithm~\ref{alg:sim}.

\begin{figure}[!t]
  \centering
  \includegraphics[width=0.5\columnwidth]{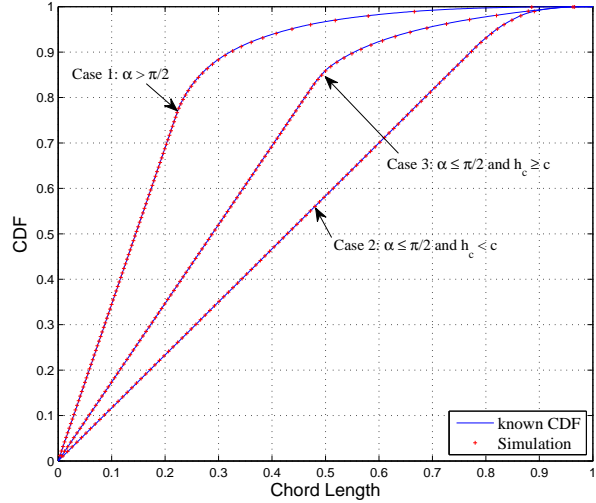}
  \caption{CLD in (\ref{eq:cld}) and simulation results for an arbitrary triangle.}
  \label{fig:cld_verify}
\end{figure}

We use three triangles listed below as examples, which cover the three cases, respectively, i.e., 1) $\alpha>\frac{\pi}{2}$: ($\alpha=\frac{130\pi}{180}, \beta=\frac{30\pi}{180}, \gamma=\frac{20\pi}{180}, a=1$); 2) $\alpha\leq\frac{\pi}{2}$ and $h_c<c$: ($\alpha=\frac{65\pi}{180}, \beta=\frac{60\pi}{180}, \gamma=\frac{55\pi}{180}, a=1$); and 3) $\alpha\leq\frac{\pi}{2}$ and $h_c\geq c$: ($\alpha=\frac{80\pi}{180}, \beta=\frac{70\pi}{180}, \gamma=\frac{30\pi}{180}, a=1$). As shown in \figurename~\ref{fig:cld_verify}, the results from the CDF in (\ref{eq:cld}) match very closely with the simulation results, verifying the correctness of (\ref{eq:cld}).

\begin{figure}[!t]
\centering
\includegraphics[width=2.5in]{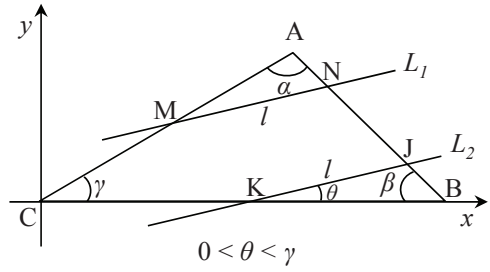}
\caption{Case example for the CLD derivation of an arbitrary triangle.}
\label{fig:cld_derive}
\end{figure}
The above CLD verification for arbitrary triangles by simulation indicates a different approach to the CLD derivation based on elementary geometry from the one in \cite{CLD}. Specifically, take the case $0<\theta<\gamma$ for example. As shown in \figurename~\ref{fig:cld_derive}, for a given chord length $l$ of the triangle, there are two corresponding lines $L_1$ and $L_2$ parallel with each other. Then we have $F(l)=\text{Pr}\{\mathcal{L}\leq l\}=\frac{||\triangle AMN||+||\triangle BJK||}{||\triangle ABC||}$, where $||\triangle||$ represents the area of a triangle and $\mathcal{L}$ the chord length random variable. The same method is applied as $\theta$ falls in other cases shown in \figurename~\ref{fig:cld_sim}. With integral over $\theta$, we can get the CLD for an arbitrary triangle. This method can easily be extended to the CLD derivation for an arbitrary convex 4-gon.
\subsubsection{Random distances distribution verification}
\begin{figure}[!t]
  \centering
  \includegraphics[width=0.5\columnwidth]{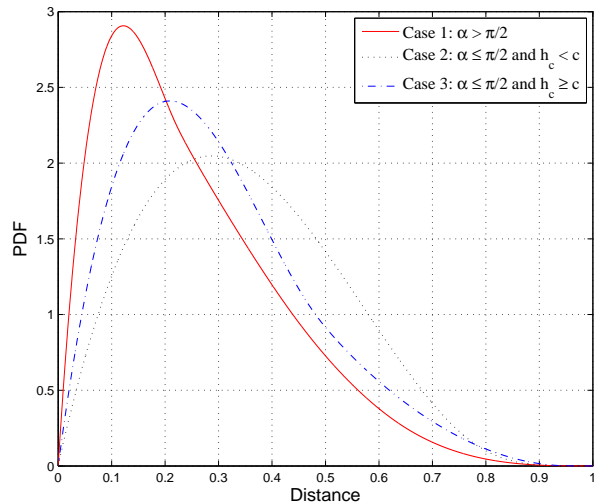}
  \caption{PDF of random distances within an arbitrary triangle.}
  \label{fig:pdd_pdf}
\end{figure}
\begin{figure}[!t]
  \centering
  \includegraphics[width=0.5\columnwidth]{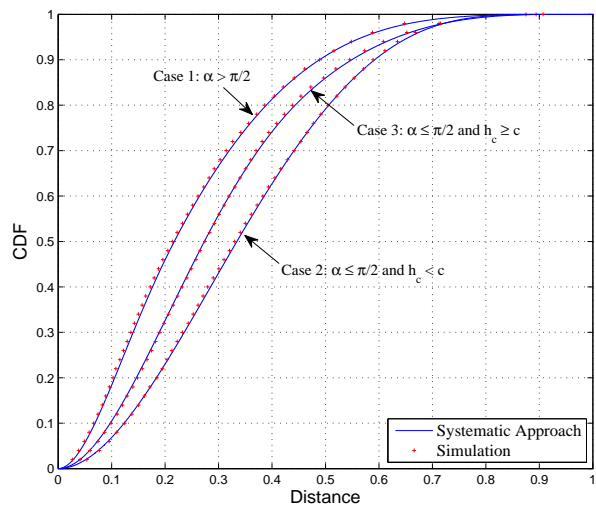}
  \caption{CDF from analysis and simulation results of random distances within an arbitrary triangle.}
  \label{fig:pdd_cdf_verify}
\end{figure}
\begin{figure}[!t]
  \centering
  \includegraphics[width=0.5\columnwidth]{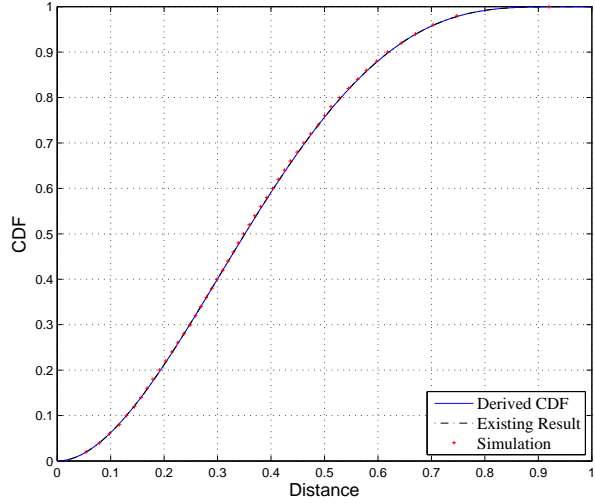}
  \caption{CDF from analysis, existing, and simulation results of random distances within an equilateral triangle.}
  \label{fig:pdd_cdf_equilateral}
\end{figure}

With the same triangle examples used above, \figurename~\ref{fig:pdd_pdf} plots the PDFs of the random distances given in (\ref{eq:g}). \figurename~\ref{fig:pdd_cdf_verify} shows a comparison between the CDFs of the random distances in (\ref{eq:G}) and the simulation results by generating 10,000 pairs of random points within the triangle with the corresponding geometric locations. It demonstrates that our distance distribution functions are very accurate when compared with the simulation results.

We also compare our derived results in (\ref{eq:G}) with existing results for equilateral triangles. In \cite{yanyan2012equilateral}, following a different approach, the PDF and CDF of random distances between two uniformly distributed points within an equilateral triangle with side length equal to 1 (which is defined as a unit equilateral triangle) are $g_{ET}(d)$ and $G_{ET}(d)$, respectively, as follows,
\begin{equation}\label{eq:pdf-equilateral}
   g_{ET}(d)= 4d\left\{
     \begin{array}{lcr}

\left(2+\frac{4\sqrt{3}\pi}{9}\right)d^2-8d+\frac{2\sqrt{3}\pi}{3}
& \quad & 0\leq d\leq \frac{\sqrt{3}}{2}\\

\frac{2\sqrt{3}}{3}\left(4d^2+6\right)\arcsin\frac{\sqrt{3}}{2d}+\left(2-
\frac{8\sqrt{3}\pi}{9}\right)d^2 \\
~~~+6\sqrt{4d^2-3}-8d-\frac{4\sqrt{3}\pi}{3}  & \quad & \frac{\sqrt{3}}{2}\leq d \leq 1\\

       0 & \quad & {\rm otherwise}
     \end{array}
   \right.,
\end{equation}
\begin{equation}\label{eq:cdf-equilateral}
   G_{ET}(d)= 2\left\{
     \begin{array}{lcl}
0 & \quad & d \leq 0 \\

\left(1+\frac{2\sqrt{3}\pi}{9}\right)d^4-\frac{16}{3}d^3+\frac{2\sqrt{3}\pi}{3}
d^2 & \quad & 0\leq d\leq \frac{\sqrt{3}}{2} \\

\frac{4\sqrt{3}d^2}{3}\left(d^2+3\right)\arcsin\frac{\sqrt{3}}{2d}+\left(\frac{
26d^2}{3}+1\right)\sqrt{d^2-\frac{3}{4}}\\
~~~+\left(1-\frac{4\sqrt{3}\pi}{9}\right)d^4-\frac{16}{3}d^3-\frac{4\sqrt{3}\pi}{3}d^2  & \quad & \frac{\sqrt{3}}{2}\leq
d \leq 1\\

       \frac{1}{2} & \quad & d \geq 1
     \end{array}
   \right..
\end{equation}%
Note that although a unit equilateral triangle with side length equal to 1 is assumed in (\ref{eq:pdf-equilateral}) and (\ref{eq:cdf-equilateral}), they can be easily scaled by a nonzero scalar, for an equilateral triangle with arbitrary side length. Specifically, let its side length be $s>0$, then
\begin{equation}\label{eq:cdf-scale}
G_{sET}(d)=P(sD\leq d)=P(D\leq \frac{d}{s})=G_{ET}(\frac{d}{s}).
\end{equation}
Therefore,
\begin{equation}\label{eq:pdf-scale}
 g_{sET}(d)=G'_{ET}(\frac{d}{s})=\frac{1}{s}g_{ET}(\frac{d}{s}).
\end{equation}
Given a unit equilateral triangle, the simplification of (\ref{eq:g}) and (\ref{eq:G}) we derived comes to the same expressions as (\ref{eq:pdf-equilateral}) and (\ref{eq:cdf-equilateral}), respectively. Meanwhile, as shown in \figurename~\ref{fig:pdd_cdf_equilateral}, the derived result for equilateral triangles in this paper matches very closely with the existing result given in (\ref{eq:cdf-equilateral}), as well as the simulation result.
\section{Random Distances between Two Arbitrary Triangles}
In this section, we first introduce the approach to obtaining the random distance distribution between two arbitrary triangles sharing a side. The results and verification for two special cases are then provided.
\subsection{Decomposition and recursion approach}
There are two cases to be discussed: 1) two arbitrary triangles forming a convex 4-gon; and 2) two arbitrary triangles forming a concave 4-gon. A decomposition and recursion approach can be applied.
\subsubsection{Two arbitrary triangles forming a convex 4-gon}
\begin{figure}[!t]
  \centering
  \includegraphics[width=2.in]{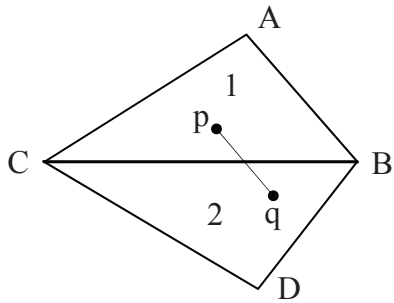}
  \caption{Random distances between two arbitrary triangles ($\triangle CAB$ and $\triangle CDB$) forming a convex 4-gon $\Box CABD$.}
  \label{fig:arbitrary-convex}
\end{figure}
As shown in \figurename~\ref{fig:arbitrary-convex}, $\triangle CAB$ labeled by 1 and $\triangle CDB$ labeled by 2 form a convex 4-gon $\Box CABD$. Let $G$, $G_1$, and $G_2$ denote the distributions of the random distances between two uniformly distributed points within $\Box CABD$, $\triangle CAB$, and $\triangle CDB$, respectively, and let $G_{12}$ denote the distribution of the random distances between two uniformly distributed points separately inside the two triangles. We have
\begin{align*}
\begin{aligned}
  G	= {} & \frac{S_1}{S}\left(\frac{S_1}{S}G_1+\frac{S_2}{S}G_{12}\right)+\frac{S_2}{S}\left(\frac{S_1}{S}G_{12}+\frac{S_2}{S}G_2\right)\,,
\end{aligned}
\end{align*}%
where $S$, $S_1$, and $S_2$ are the areas of $\Box CABD$, $\triangle CAB$, and $\triangle CDB$, respectively. Therefore,
\begin{align}\label{eq:arbitrary-convex}
\begin{aligned}
  G_{12} = {} & \frac{S^2G-S_1^2G_1-S_2^2G_2}{2S_1S_2}\,.
\end{aligned}
\end{align}
$G_1$ and $G_2$ have been given in (\ref{eq:G}). With the CLD of an arbitrary convex 4-gon, which can be derived using the approach introduced in Section \ref{subsec:verify}, $G$ can be obtained by utilizing (\ref{eq:G}), and the author in \cite{CLD_PARALLEL} discussed the CLD for parallelograms. As a case study, $G_{12}$ will be derived in Section \ref{sec:2T_convex_verify} when $\Box CABD$ is a rhombus formed by two congruent isosceles triangles, where $G$ for a rhombus is known in \cite{yanyan2011rhombus}.
\subsubsection{Two arbitrary triangles forming a concave 4-gon}
\begin{figure}[!t]
  \centering
  \includegraphics[width=2.0in]{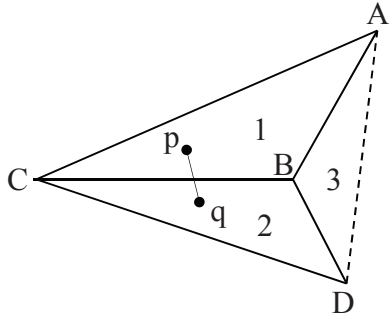}
  \caption{Random distances between two arbitrary triangles ($\triangle CAB$ and $\triangle CDB$) forming a concave 4-gon $\Box CABD$.}
  \label{fig:arbitrary-concave}
\end{figure}
If $\Box CABD$ is a concave 4-gon, as shown in \figurename~\ref{fig:arbitrary-concave}, (\ref{eq:arbitrary-convex}) still applies. The approach based on CLD to obtaining $G$ may not apply to the concave 4-gon due to the ambiguity in its chord length definition and in the relationship between its CLD and point distance distribution. 

By linking $AD$ with the triangle $\triangle ABD$ labeled by 3 as shown in \figurename~\ref{fig:arbitrary-concave}, the decomposition and recursion approach can still apply, and we have
\begin{align*}
\begin{aligned}
  G	= {} & \frac{S_1}{S}\left(\frac{S_1}{S}G_1+\frac{S_2}{S}G_{12}+\frac{S_3}{S}G_{13}\right)+\frac{S_2}{S}\left(\frac{S_2}{S}G_{2}+\frac{S_1}{S}G_{12}+\frac{S_3}{S}G_{23}\right)\\
  {}&+\frac{S_3}{S}\left(\frac{S_3}{S}G_3+\frac{S_1}{S}G_{13}+\frac{S_2}{S}G_{23}\right)\,,
\end{aligned}
\end{align*}
where $G$ is the distribution of the random distances between two uniformly distributed points within $\triangle ACD$. If $\triangle ACD$ is an equilateral triangle, and triangle 1, 2, and 3 are congruent with each other, $G_1=G_2=G_3$, $G_{12}=G_{13}=G_{23}$, $S_1=S_2=S_3=\frac{1}{3}S$, and thus $G_{12}=\frac{1}{2}(3G-G_1)$, where the $G$ of an equilateral triangle is given in (\ref{eq:cdf-equilateral}). The detailed results will be given in Section \ref{sec:2T_concave_verify}.

\subsection{Results and verification}
In this subsection, we provide two case studies. One is the distance distribution between two isosceles triangles forming a rhombus, and the other is that when two obtuse isosceles triangles with the acute angle equal to $\frac{\pi}{6}$  form a concave 4-gon. For the first case, when the two triangles are equilateral ones, the results have been given in \cite{yanyan2012equilateral}, and thus we only focus on the case when the two triangles are obtuse ones with the acute angle equal to $\frac{\pi}{6}$.
\subsubsection{Distance distribution between two isosceles triangles (with the acute angle equal to $\frac{\pi}{6}$) forming a rhombus}\label{sec:2T_convex_verify}
\begin{figure}[!t]
  \centering
  \includegraphics[width=2.5in]{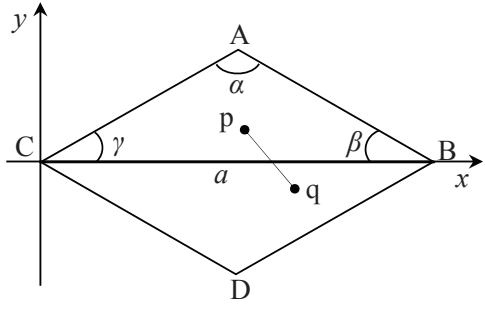}
  \caption{Random distances between two unit congruent triangles ($\triangle CAB$ and $\triangle CDB$) forming a rhombus $\diamond CABD$.}
  \label{fig:rhombus}
\end{figure}
According to \cite{yanyan2011rhombus}, the PDF and CDF of random distances between two uniformly distributed points within a unit rhombus are $g_R(d)$ and $G_R(d)$, respectively, as follows,
\begin{equation}\label{eq:fd_r_within}
  g_R(d)=2d\left\{
    \begin{array}{lr}

\left(\frac{4}{3}+\frac{2\sqrt{3}\pi}{27}\right)d^2-\frac{16}{3}d+\frac{2\sqrt{3}\pi}{3} & 0\leq d\leq \frac{\sqrt{3}}{2}\\

\frac{8\sqrt{3}}{3}\left(1+\frac{d^2}{3}\right)\arcsin\frac{\sqrt{3}}{2d}
+\left(\frac{4}{3}-\frac{10\sqrt{3}\pi}{27}\right)d^2-\frac{16}{3}d+\frac{10}{3}
\sqrt{4d^2-3}-\frac{2\sqrt{3}\pi}{3} & \frac{\sqrt{3}}{2}\leq d\leq 1\\

\frac{4\sqrt{3}}{3}\left(1-\frac{d^2}{3}\right)\arcsin\frac{\sqrt{3}}{2d}
-\left(\frac{2}{3}-\frac{2\sqrt{3}\pi}{27}\right)d^2+\sqrt{4d^2-3}-\frac{2\sqrt{3}\pi}{9}-1 & 1\leq d\leq \sqrt{3} \\

      0 & {\rm otherwise}
    \end{array}
  \right.,
\end{equation}
\begin{equation}\label{eq:Fd_r_within}
  G_R(d)=\left\{
    \begin{array}{lr}

\left(\frac{2}{3}+\frac{\sqrt{3}\pi}{27}\right)d^4-\frac{32}{9}d^3+\frac{
2\sqrt{3}\pi}{3}d^2 & 0\leq d\leq \frac{\sqrt{3}}{2}\\

\frac{4\sqrt{3}}{3}\left(2d^2+\frac{d^4}{3}\right)\arcsin\frac{\sqrt{3}}{2d}+\left(\frac{2}{3}-\frac{5\sqrt{3}\pi}{27}\right)d^4-\frac{32
}{9}d^3-\frac{2\sqrt{3}\pi}{3}d^2 \\
~~~~+\frac{1}{6}(14d^2+3)\sqrt{4d^2-3} & \frac{\sqrt{3}}{2}\leq d\leq 1\\

\frac{2\sqrt{3}}{3}\left(2d^2-\frac{d^4}{3}\right)\arcsin\frac{\sqrt{3}}{
2d}+\left(\frac{\sqrt{3}\pi}{27}-\frac{1}{3}\right)d^4-\left(\frac{2\sqrt{3}\pi}{9}+1\right)d^2 \\
~~~~+\frac{1}{36}(22d^2+15)\sqrt{4d^2-3}+\frac{1}{4} & 1\leq d\leq \sqrt{3} \\

1 & d\geq \sqrt{3}
    \end{array}
 \right..
\end{equation}

Given a unit isosceles triangle $T: (\alpha=\frac{120\pi}{180}, \beta=\frac{30\pi}{180}, \gamma=\frac{30\pi}{180}, a=1)$, we can derive the random distances between two $T$s forming a rhombus with side length equal to $\frac{\sqrt{3}}{3}$, such as $|pq|$ shown in \figurename~\ref{fig:rhombus}. Denote the PDF and CDF of $|pq|$ by $g^r_{2T}(d)$ and $G^r_{2T}(d)$, respectively. Following \cite{yanyan2012equilateral}, we know that
\begin{equation*}
\left\{
  \begin{array}{rl}
  \sqrt{3}g_R(\sqrt{3}d)=&\frac{1}{2}g_T(d)+\frac{1}{2}g^r_{2T}(d)\\
  G_R(\sqrt{3}d)=&\frac{1}{2}G_T(d)+\frac{1}{2}G^r_{2T}(d)
  \end{array}
  \right.,
\end{equation*}
and thus we have
\begin{equation*}
\left\{
  \begin{array}{rl}
  g^r_{2T}(d)=&2\sqrt{3}g_R(\sqrt{3}d)-g_T(d)\\
  G^r_{2T}(d)=&2G_R(\sqrt{3}d)-G_T(d)
  \end{array}
  \right.,
\end{equation*}
where $g_T(d)$ and $G_T(d)$ have been given in (\ref{eq:g}) and (\ref{eq:G}), respectively. Therefore,
\begin{equation}\label{eq:pdf-rT2}
  g^r_{2T}(d)=\frac{8}{3}d\left\{
    \begin{array}{lr}
    36d-(9+13\sqrt{3}\pi)d^2 & 0\leq d\leq\frac{\sqrt{3}}{6} \\

    \left(36\sqrt{3}\arccos\frac{\sqrt{3}}{6d}-13\sqrt{3}\pi-9 \right)d^2-9\sqrt{36d^2-3}\\
    ~~~-6\sqrt{3}\arcsin\frac{\sqrt{3}}{6d}+3\sqrt{3}\pi+36d& \frac{\sqrt{3}}{6}\leq d\leq\frac{1}{2} \\

    12\sqrt{3}(1+d^2)\arcsin\frac{1}{2d}+ (36\sqrt{3}\arccos\frac{\sqrt{3}}{6d}-9-19\sqrt{3}\pi) d^2\\
     ~~~+15\sqrt{12d^2-3}-3\sqrt{3}\pi-6\sqrt{3}\arcsin\frac{\sqrt{3}}{6d}\\
     ~~~-9\sqrt{36d^2-3}+36d& \frac{1}{2}\leq d\leq\frac{\sqrt{3}}{3} \\

    (12\sqrt{3}\arccos\frac{1}{2d}-6\sqrt{3}\arcsin{\frac{1}{2d}}-3\sqrt{3}\pi -18)d^2\\
    ~~~-\frac{9\sqrt{3}}{2}\sqrt{4d^2-1}+36d-\frac{9}{2} & \frac{\sqrt{3}}{3}\leq d\leq 1 \\
    0 & {\rm otherwise}
    \end{array}
  \right.,
\end{equation}
\begin{equation}\label{eq:cdf-rT2}
  G^r_{2T}(d)=\left\{
    \begin{array}{lr}
     32d^3-\left(\frac{26\sqrt{3}\pi}{3}+6\right)d^4& 0\leq d\leq\frac{\sqrt{3}}{6} \\

     \left(24\sqrt{3}\arccos\frac{\sqrt{3}}{6d}-\frac{26\sqrt{3}\pi}{3}-6 \right)d^4 +32d^3 \\
      ~~~+\left( 4\sqrt{3}\pi-\frac{26}{3}\sqrt{36d^2-3}-8\sqrt{3}\arcsin\frac{\sqrt{3}}{6d} \right)d^2-\frac{1}{9}\sqrt{36d^2-3}& \frac{\sqrt{3}}{6}\leq d\leq\frac{1}{2} \\

     \frac{1}{9}(126d^2+9)\sqrt{12d^2-3}-\frac{1}{9}(1+78d^2)\sqrt{36d^2-3}\\
     ~~~-\frac{38d^2}{3}\left(\frac{12\sqrt{3}}{19}\arcsin\frac{\sqrt{3}}{6d}-\frac{36\sqrt{3}d^2}{19}\arccos\frac{\sqrt{3}}{6d}\right.\\
     ~~~\left.-\frac{12\sqrt{3}}{19}(2+d^2)\arcsin\frac{1}{2d}+\sqrt{3}\pi\left(d^2+\frac{6}{19}\right)+\frac{9d^2}{19}-\frac{48d}{19}\right)& \frac{1}{2}\leq d\leq\frac{\sqrt{3}}{3} \\

     \left(8\sqrt{3}\arccos\frac{1}{2d}-4\sqrt{3}\arcsin\frac{1}{2d}-2\sqrt{3}\pi-12\right)d^4\\
     ~~~+32d^3+\frac{\sqrt{3}}{2}(1-10d^2)\sqrt{4d^2-1} -6d^2+\frac{1}{2} & \frac{\sqrt{3}}{3}\leq d\leq 1 \\
    1 & {\rm d\geq 1}
    \end{array}
  \right..
\end{equation}
\begin{figure}[!td]
  \centering
  \includegraphics[width=0.5\columnwidth]{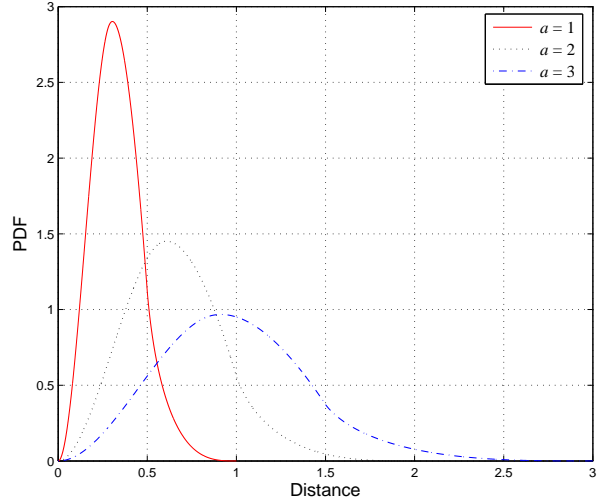}
  \caption{PDF of random distances between two adjacent triangles forming a rhombus.}
  \label{fig:2T-Rhombus-pdf}
\end{figure}
\begin{figure}[!t]
  \centering
  \includegraphics[width=0.5\columnwidth]{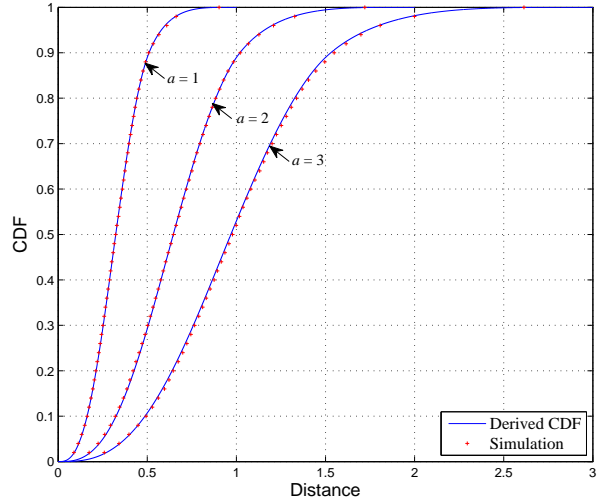}
  \caption{CDF from analysis and simulation results of random distances between two adjacent triangles forming a rhombus.}
  \label{fig:2T-Rhombus-cdf}
\end{figure}

\figurename~\ref{fig:2T-Rhombus-pdf} plots the PDF given in (\ref{eq:pdf-rT2}) for $a=1,2,~\text{and}~3$ (according to (\ref{eq:pdf-scale}), $g^r_{a2T}(d)=\frac{1}{a}g^r_{2T}(\frac{d}{a})$). \figurename~\ref{fig:2T-Rhombus-cdf} shows a comparison between our derived CDF given in (\ref{eq:cdf-rT2}) and the simulation results by generating 10,000 pairs of random pints, with the two points of each pair falling in the two triangles, respectively, such as $p$ and $q$ shown in \figurename~\ref{fig:rhombus}. It demonstrates that our derived distance distribution functions are very accurate when compared with the simulation results.

\subsubsection{Distance distribution between two isosceles triangles (with the acute angle equal to $\frac{\pi}{6}$) forming a concave 4-gon}\label{sec:2T_concave_verify}
\begin{figure}[!t]
  \centering
  \includegraphics[width=2.5in]{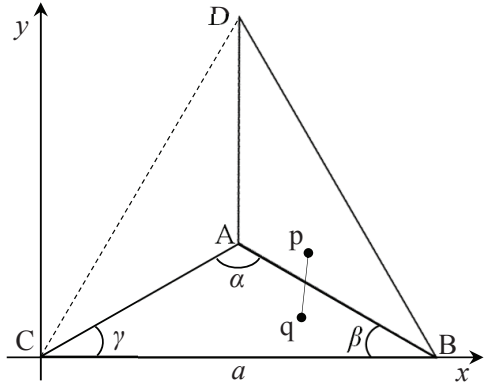}
  \caption{Random distances between two unit congruent triangles ($\triangle CAB$ and $\triangle ADB$) forming a concave 4-gon $\Box CADB$.}
  \label{fig:concave}
\end{figure}
The two unit isosceles triangle $T: (\alpha=\frac{120\pi}{180}, \beta=\frac{30\pi}{180}, \gamma=\frac{30\pi}{180}, a=1)$ can also form a concave 4-gon, such as $\Box CADB$ shown in \figurename~\ref{fig:concave}. By linking $CD$, the PDF $g^c_{2T}(d)$ and CDF $G^c_{2T}(d)$ of the random distances between the two $T$s, such as $|pq|$ shown in \figurename~\ref{fig:concave}, can be obtained using the following equations,
\begin{equation*}
\left\{
  \begin{array}{rl}
  g_{ET}(d)=&\frac{1}{3}g_T(d)+\frac{2}{3}g^c_{2T}(d)\\
  G_{ET}(d)=&\frac{1}{3}G_T(d)+\frac{2}{3}G^c_{2T}(d)
  \end{array}
  \right.,
\end{equation*}
and thus
\begin{equation*}
\left\{
  \begin{array}{rl}
  g^c_{2T}(d)=& \frac{1}{2}\left(3g_{ET}(d)-g_T(d)\right)\\
  G^c_{2T}(d)=&\frac{1}{2}\left(3G_{ET}(d)-G_T(d)\right)
  \end{array}
  \right.,
\end{equation*}
where $g_T(d)$ and $G_T(d)$ have been given in (\ref{eq:g}) and (\ref{eq:G}), respectively; $g_{ET}(d)$ and $G_{ET}(d)$ are in (\ref{eq:pdf-equilateral}) and (\ref{eq:cdf-equilateral}), respectively. Therefore,
\begin{equation}\label{eq:pdf-cT2}
  g^c_{2T}(d)=\left\{
    \begin{array}{lr}
    \frac{8}{3}(12\sqrt{3}-9d-6\sqrt{3}\pi d)d^2& 0\leq d\leq\frac{\sqrt{3}}{6} \\

    \frac{4}{3}d\left[\left(36\sqrt{3}\arccos\frac{\sqrt{3}}{6d}-12\sqrt{3}\pi-18\right)d^2+24\sqrt{3}d\right.\\
    ~~~\left.-9\sqrt{36d^2-3}+3\sqrt{3}\pi-6\sqrt{3}\arcsin\frac{\sqrt{3}}{6d} \right]& \frac{\sqrt{3}}{6}\leq d\leq\frac{\sqrt{3}}{3} \\

    \frac{4}{9}d\left[ \left(36\sqrt{3}\arccos\frac{1}{2d}-6\sqrt{3}\pi\right)d^2-27\sqrt{12d^2-3}\right.\\
    ~~~\left.-18\sqrt{3}\arcsin\frac{1}{2d}+12\sqrt{3}\pi \right]& \frac{\sqrt{3}}{3}\leq d\leq\frac{\sqrt{3}}{2} \\

    \frac{32\sqrt{3}}{3}d\left[\frac{9\sqrt{3}}{8}\sqrt{4d^2-3}-\frac{9}{8}\sqrt{4d^2-1}+\left(\frac{3}{2}d^2+\frac{9}{4}\right)\arcsin\frac{\sqrt{3}}{2d}\right.\\
    ~~~\left.+\frac{3}{2}d^2\arccos\frac{1}{2d}-\pi d^2-\frac{3}{4}\arcsin\frac{1}{2d}-\frac{5}{8}\pi  \right]& \frac{\sqrt{3}}{2}\leq d\leq 1 \\

    0 & {\rm otherwise}
    \end{array}
  \right.,
\end{equation}%
\begin{equation}\label{eq:cdf-cT2}
  G^c_{2T}(d)=\left\{
    \begin{array}{lr}
    \frac{32\sqrt{3}}{3}d^3-(6+4\sqrt{3}\pi)d^4& 0\leq d\leq\frac{\sqrt{3}}{6} \\

    \frac{1}{18}(-78d^2-1)\sqrt{36d^2-3}-4d^2\left[ \sqrt{3}\arcsin\frac{\sqrt{3}}{6d}\right.\\
    ~~~\left.-3\sqrt{3}d^2\arccos\frac{\sqrt{3}}{6d}+\left(\pi d^2 -\frac{8}{3}d-\frac{\pi}{2}\right)\sqrt{3}+\frac{3}{2}d^2 \right]& \frac{\sqrt{3}}{6}\leq d\leq\frac{\sqrt{3}}{3} \\

    \frac{\sqrt{3}}{6}\left[16\pi d^2-4\pi d^4-\sqrt{4d^2-1}(1+26d^2)\right.\\
    ~~~\left.-24d^2\left(\arcsin\frac{1}{2d}-d^2\arccos\frac{1}{2d}\right) \right]& \frac{\sqrt{3}}{3}\leq d\leq\frac{\sqrt{3}}{2} \\

    \frac{1}{6}(78d^2+9)\sqrt{4d^2-3}-\frac{8\sqrt{3}}{3}\left[\left(\frac{1}{16}+\frac{13}{8}d^2\right)\sqrt{4d^2-1}\right.\\
    ~~~\left.+d^2\left(\frac{3}{2}\arcsin\frac{1}{2d}-\left(\frac{3}{2}d^2+\frac{9}{2}\right)\arcsin\frac{\sqrt{3}}{2d}\right.\right.\\
    ~~~\left.\left.-\frac{3}{2}d^2\arccos\frac{1}{2d}+\pi d^2+\frac{5}{4}\pi \right) \right] & \frac{\sqrt{3}}{2}\leq d\leq 1 \\

    1 & {\rm d\geq 1}
    \end{array}
  \right..
\end{equation}%
\begin{figure}[!t]
  \centering
  \includegraphics[width=0.5\columnwidth]{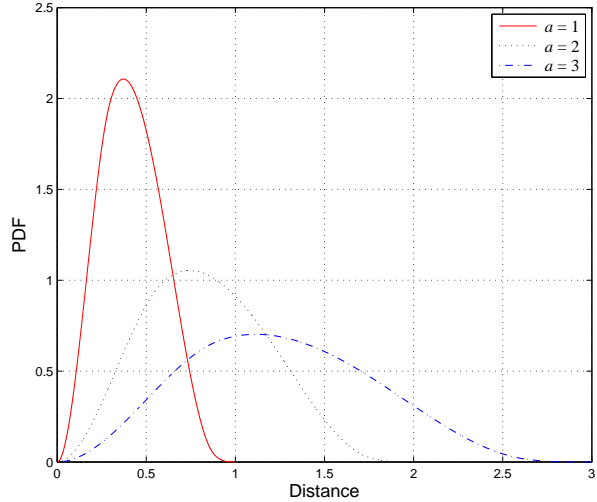}
  \caption{PDF of random distances between two adjacent triangles forming a concave 4-gon.}
  \label{fig:2T-Concave-pdf}
\end{figure}%
\begin{figure}[!t]
  \centering
  \includegraphics[width=0.5\columnwidth]{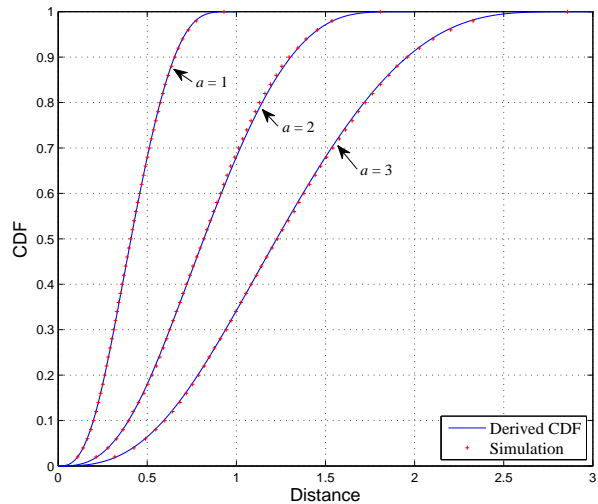}
  \caption{CDF from analysis and simulation results of random distances between two adjacent triangles forming a concave 4-gon.}
  \label{fig:2T-Concave-cdf}
\end{figure}%
The PDFs for cases when $a=1,2,~\text{and}~3$ are plotted in \figurename~\ref{fig:2T-Concave-pdf}. The corresponding CDFs are compared with the simulation results, as shown in \figurename~\ref{fig:2T-Concave-cdf}, which demonstrates that our derived results are very accurate.
\section{Discussions and Conclusions}\label{sec:discussion}
In this work, we proposed a systematic approach to obtaining the random distance distributions between two uniformly distributed points associated with arbitrary triangles. The two points can either be inside the same triangle or separately inside two adjacent triangles sharing a side. For the former, the results are obtained based on the CLD for arbitrary triangles. For the latter, a systematic approach based on decomposition and recursion is applied, by separately discussing the formed convex and concave 4-gon. The results for two special cases, i.e., two congruent isosceles triangles with the acute angle equal to $\frac{\pi}{6}$ forming a rhombus and a concave 4-gon, are given in detail.

\section*{Acknowledgment}
This work is supported in part by the NSERC, CFI and BCKDF.

\end{document}